\documentclass[11pt]{article}
\usepackage{amsmath, amsthm, amssymb}
\usepackage[dvips]{graphics,color}
\usepackage{epsfig}
\usepackage{latexsym}
\usepackage{graphicx}

\begin{document}

\swapnumbers
\theoremstyle{plain}% default
\newtheorem{thm}{Theorem}[section]
\newtheorem{lem}[thm]{Lemma}
\newtheorem{prop}[thm]{Proposition}
\newtheorem{cor}[thm]{Corollary}

\theoremstyle{definition}
\newtheorem{defn}[thm]{Definition}
\newtheorem{conj}[thm]{Conjecture}
\newtheorem{exmp}[thm]{Example}
\newtheorem{nota}[thm]{Notation}

\theoremstyle{remark}
\newtheorem*{rem}{Remark}
\newtheorem*{note}{Note}

\title{Unit Rectangle Visibility Graphs}
\author{
   \emph{Alice M. Dean}\\Department of Mathematics and Computer Science\\
            Skidmore College, Saratoga Springs, NY, adean@skidmore.edu
\\\\
    \emph{Joanna A. Ellis-Monaghan}\footnote{Research supported by VT EPSCoR under grant NSF EPS 0236876 
and the National Security Agency.}
    \\Department of Mathematics\\
        Saint Michael's College, Colchester, VT, jellis-monaghan@smcvt.edu
\\\\
    \emph{Sarah Hamilton}\footnote{Research supported by
            NASA under Training Grant NGT5-40110 to the Vermont Space Grant
            Consortium.}
    \\Department of Mathematics\\
            Saint Michael's College, Colchester, VT, shamilton2@smcvt.edu
\\\\
    \emph{Greta Pangborn}\footnote{Research supported by VT EPSCoR under grant NSF EPS 0236876
    and by NASA under Training Grant NGT5-40110 to the Vermont Space Grant
            Consortium.}
    \\Department of Computer Science \\
            Saint Michael's College, Colchester, VT, gpangborn@smcvt.edu
}

%\date{DRAFT: \noindent Current Version for Authors' Eyes Only
%            \\Revised: 11 August 2007}

\maketitle

\begin{abstract}
Over the past twenty years, rectangle visibility graphs have generated
considerable interest, in part due to their applicability to VLSI chip design.
Here we study unit rectangle visibility graphs, with fixed dimension
restrictions more closely modeling the constrained
dimensions of gates and other circuit components in
computer chip applications.
A graph $G$ is a unit rectangle visibility graph
(URVG) if its vertices can be represented by closed unit squares in
the plane with sides parallel to the axes and pairwise disjoint interiors,
in such a way that two vertices are adjacent if and
only if there is a non-degenerate horizontal or vertical band of visibility joining the two
rectangles. Our results include necessary and sufficient conditions for
$K_n$, $K_{m,n}$, and trees to be URVGs, as well as a number of general edge bounds.
\end{abstract}

\section{Introduction}\label{sec:intro}

Over the past twenty years the difficulty of VLSI
chip design and layout problems has motivated the study of bar visibility
graphs (BVGs) \cite{DGH04, DV03, DHLM83, RT86, SLMW85, TT86, Wismath85} and their
two-dimensional counterparts, rectangle visibility graphs (RVGs)
\cite{BDHS97, DH94, HSV99, Shermer96, Shermer96b, SW03}. In these constructions, horizontal
bars or rectangles in the plane model gates or other chip components, and edges are
modeled by vertical visibilities between bars, or by vertical and horizontal visibilities
between rectangles. The two visibility directions in  RVGs provide a model for
two-layer chips with wires running horizontally on one layer and vertically on the
other.  The dimensions of bars and rectangles in BVGs and RVGs may vary arbitrarily,
but chip components typically have restricted area and aspect ratios. In order to
more closely model the restricted dimensions of chip components, Dean and Veytsel
\cite{DV03} studied unit bar visibility graphs (UBVGs), in which all bars have
equal length. A related model, using boxes in 3-space, was studied in
\cite{BJMO94, FM99}.  In this paper we study the similarly restricted class of RVGs,
unit rectangle visibility graphs (URVGs), in which all rectangles are unit
squares.

In Section~\ref{sec:defns} we give definitions and basic properties that we
use throughout the paper. In Section~\ref{sec:complete} we characterize the
complete graphs that are URVGs. In Section~\ref{sec:trees} we characterize URVG
trees, and we show that any graph with linear arboricity 2 is a URVG.
In section~\ref{sec:bipartite} we characterize which complete
bipartite graphs are URVGs or subgraphs of URVGs.
In section~\ref{sec:bounds} we give edge bound results for URVGs
as well as examples that show these bounds are
tight up to constant coefficients. We provide edge bounds for depth-$s$ UBV and
URV trees, bipartite URVGs, and arbitrary URVGs. In Section~\ref{sec:concl} we
 state two open problems on URVGs.

\section{Definitions and Basic Properties} \label{sec:defns}

 \begin{defn} A graph $G$ is a {\em unit rectangle visibility graph} or
  {\em URVG} if its vertices can be represented by closed unit squares in
  the plane with sides parallel to the axes and pairwise disjoint
  interiors, in such a way that two vertices are adjacent if and
  only if there is an unobstructed non-degenerate (positive width)
  horizontal or vertical band
  of visibility joining the two rectangles.
\end{defn}

We denote the square in the URV layout corresponding to a vertex $v$ by $S_v$.
We identify the position of the square $S_v$ in a
URV layout by its bottom-left corner coordinates $(x_v,y_v)$.
We define $X_v$ to be the line segment given by the intersection of
the line $x=x_v$ with the square $S_v$, and $Y_v$ to be the line segment
given by the intersection of the line $y=y_v$ with the square $S_v$.

\begin{figure}[htb]%[htbp]
  \centerline{\epsfig{file=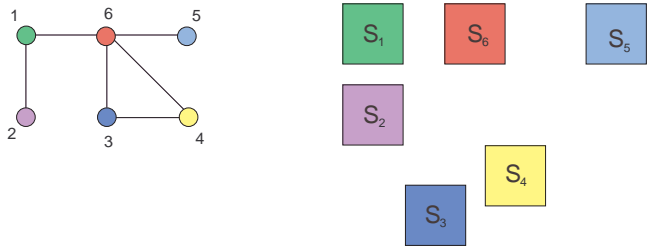, width=0.7\textwidth}}
  \caption{A  graph with a URV layout.}
  \label{fig:URVexample}
\end{figure}

Two squares $S_v$ and $S_w$ are called {\em flush} if $x_v=x_w$ or $y_v=y_w$
(this does not preclude other squares obstructing visibility between $S_v$ and $S_w$).
In Fig.~\ref{fig:URVexample}, squares $S_2$ and
$S_3$ are collinear but not flush, and squares $S_1$ and $S_2$ are flush,
as are squares $S_1$, $S_5$, and $S_6$.

\begin{defn} A graph $G$ is a {\em weak} unit rectangle visibility graph
if its vertices can be represented by closed unit squares in
  the plane with sides parallel to the axes and pairwise disjoint interiors,
  in such a way that whenever two vertices are adjacent
  there is an unobstructed non-degenerate, horizontal or vertical
  band of visibility joining the two rectangles.  Equivalently,
  $G$ is a weak URVG if it is a subgraph of a URVG.
 \end{defn}

\begin{figure}[htb]%[htbp]
  \centerline{\epsfig{file=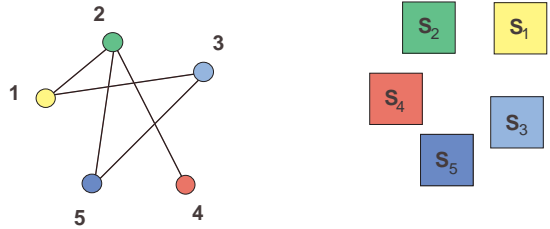, width=0.6\textwidth}}
  \caption{A graph with a weak URV layout.}
  \label{fig:URVGweak}
\end{figure}

An example of a graph with a weak URV layout is given in Fig.~\ref{fig:URVGweak}.
(There is a band of visibility between squares $S_3$ and $S_4$ but no edge $\{3,4\}$.)
We conclude this section with a straightforward but useful proposition.

\begin{prop}\label{prop:UBVGunion}
If $G$ has a URVG layout $L$, let  $G_X$ and  $G_Y$ denote the graphs
induced by the horizontal and vertical visibilities of $L$.
$G_X$ and $G_Y$ are UBVGs with bars given by $\{X_v | v \in G\}$ and 
$\{Y_v | v \in G\}$, and $G=G_X \bigcup G_Y$.
\end{prop}

\begin{figure}[htb]%[htbp]
  \centerline{\epsfig{file=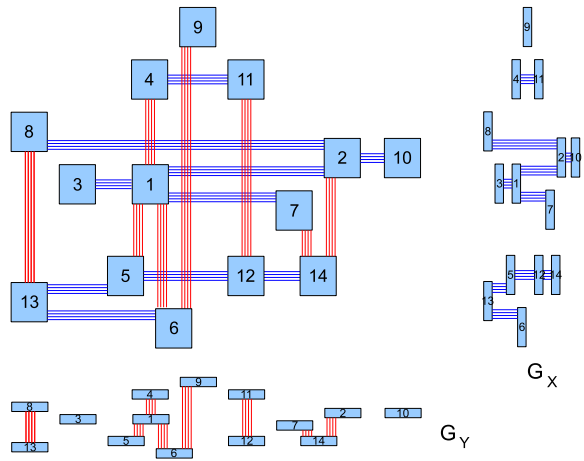, width=0.6\textwidth}}
  \caption{URV layout and corresponding UBV layouts for $G_X$ and $G_Y$.}
  \label{fig:ubvgunion}
\end{figure}

Fig.~\ref{fig:ubvgunion} illustrates the decomposition of a URV layout
into horizontal and vertical UBV layouts.

\section{Cycles and Complete Graphs}\label{sec:complete}

In this section we characterize the cycles and complete graphs that are URVGs.

 \begin{prop} \label{prop: infinite length cycle}
 The $n$-cycle $C_n$ is a URVG.
 \end{prop}

\begin{figure}[htb]%[htbp]
  \centerline{\epsfig{file=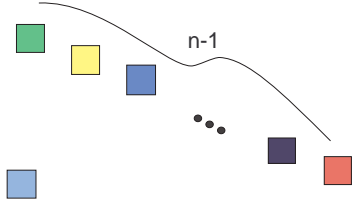, width=0.4\textwidth}}
  \caption{A layout of the $n$-cycle $C_n$.}
  \label{fig:NCcycle}
\end{figure}

See Fig.~\ref{fig:NCcycle} for a layout of $C_n$.

The following theorem of Erd\"{o}s and Szekeres \cite{ES35} will
be used in characterizing which graphs are URVGs.

 \begin{figure}[htb]%[htbp]
  \centerline{\epsfig{file=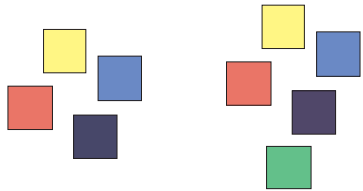, width=0.4\textwidth}}
  \caption{URV layouts of $K_4$ and of $K_5$ minus one edge.}
\label{fig:k4k5}
\end{figure}

 \begin{thm} \label{thm:subseq}(Erd\"{o}s and Szekeres \cite{ES35})
 For $n\ge0$, every sequence $a_1, a_2,$ \ldots, $a_{n^2+1}$ of $n^2+1$
    terms contains a monotonic subsequence of $n+1$ terms.
 \end{thm}

 \begin{thm} \label{thm:K5}
 $K_n$ is a URVG if and only if $n \leq 4$.
 \end{thm}

 \begin{proof}
 We first note that since all edges are present in a complete graph,
 any URV layout of $K_n$ gives a URV layout of $K_m$  for all $m \leq n$.
 Thus it suffices to prove that $K_4$ is a URVG  and $K_5$ is not a URVG.
 Fig.~\ref{fig:k4k5} gives a URV layout of $K_4$.
 To show that $K_5$ is not a URVG, we first observe that if squares
 $S_1, S_2, S_3$ give a URV layout of $K_3$ with $x_1 \leq x_2 \leq x_3$ then
 $(y_1, y_2, y_3)$ must be non-monotonic.  This follows since if the sequence
 $(y_1, y_2, y_3)$ is monotonic, then $S_2$ blocks $S_3$ from seeing $S_1$.
 Now suppose we have a URV layout of five squares.
 By relabeling if necessary, we
 may assume that $x_1 \leq x_2 \leq x_3 \leq x_4 \leq x_5$.  Now consider
 the sequence of $y$-coordinates.  By Theorem~\ref{thm:subseq} this sequence
 has a monotonic subsequence of length 3, say $(y_{i_1}, y_{i_2}, y_{i_3}).$
 Thus the squares $S_{i_1}, S_{i_2}, S_{i_3}$ cannot be a layout of $K_3$,
 and the five squares cannot be a layout of $K_5$.
\end{proof}

\begin{cor} \label{cor:K5sub}
Any graph G that contains $K_5$ as a subgraph is not a URVG.
\end{cor}

\begin{figure}[htb]%[htbp]
 \centerline{\epsfig{file=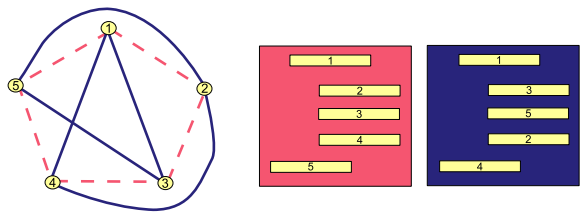, width=0.6\textwidth}}
\caption{$K_5$ is the union of two UBVGs.}
\label{fig:k5ubvgunion}
\end{figure}

\begin{rem}
Proposition \ref{prop:UBVGunion} states that every URVG is the union of two UBVGs.
Fig.~\ref{fig:k5ubvgunion} gives a decomposition of $K_5$ into the union of two
UBVGs, so the converse of Proposition~\ref{prop:UBVGunion} is false.
\end{rem}

\section{Trees, Caterpillars, and Arboricity 2}\label{sec:trees}

In this section we characterize the trees that are URVGs. We begin with
a simple necessary degree condition.

\begin{thm}\label{thm:deg 7-cycle}
If $v$ is a vertex of a URVG $G$ with degree $\ge 7$,
    then $v$ lies on a cycle.
\end{thm}

\begin{proof}
Suppose $G$ is a URVG having a vertex $v$ with degree $\ge 7$. Let $G_X$ and $G_Y$ be the UBVGs induced by the
horizontal and vertical visibilities of a URV layout of $G$, as
described in Proposition \ref{prop:UBVGunion}. Since $deg(v)\ge7$, it follows
from the Pigeonhole Principle that $deg(v)\ge4$ in $G_X$ or $G_Y$. Without loss
of generality, assume that $deg(v)\ge4$ in $G_X$. In 
(\cite{DV03}, Corollary 1.7) it is shown that
any vertex in a UBVG with degree $\ge 4$ lies on a cycle. Since $G_X$ is a
subgraph of $G$, it follows that $v$ lies on a cycle in $G$.
\end{proof}

\begin{cor} A URVG tree T has maximum degree~$\leq~6$.
    \label{cor: tree max deg 6}
\end{cor}

\begin{figure}[htb]%[htbp]
 \centerline{\epsfig{file=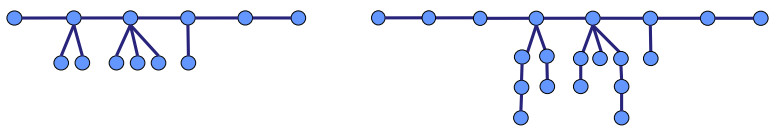, width=0.8\textwidth}}
 \caption{A caterpillar and a subdivided caterpillar.}
 \label{fig:caterpillar}
\end{figure}

\begin{defn}
A \emph{caterpillar} is a tree in which all vertices with degree greater than
1 lie on a single path. Such a path is called a \emph{spine} of the caterpillar if it has maximal length. 
A \emph{subdivided caterpillar} is a caterpillar
in which each edge may be replaced by a path of arbitrary length. 
A {\em leg} of a caterpillar or subdivided caterpillar is a path having one endpoint on the spine and the other a degree-1 vertex.
See Fig.~\ref{fig:caterpillar}.
\end{defn}

The following theorem of Dean and Veytsel \cite{DV03} characterizes the trees that are UBVGs.

\begin{thm}(Dean and Veytsel \cite{DV03}, Theorem 3.5) \label{thm:UBVGtree}
A tree is a UBVG if and only if it is a subdivided caterpillar with maximum degree 3.
\end{thm}

It follows from Proposition~\ref{prop:UBVGunion}
that if a tree is a URVG, it is the union of two subdivided caterpillar forests,
each with maximum degree 3. 
Although in general it is not true that the union of two UBVGs is a URVG, as
illustrated in Fig.~\ref{fig:k5ubvgunion}, we show that a tree $T$ is a URVG if and
only if its horizontal and vertical UBVG subgraphs, $T_X$ and $T_Y$, are UBV forests.

\begin{thm}\label{thm:URVGtree}
A tree $T$ is a URVG if and only if it is the union of two subdivided caterpillar
forests, each with maximum degree 3. 
\end{thm}

\begin{proof}
Necessity follows from Theorem~\ref{thm:UBVGtree} and Proposition~\ref{prop:UBVGunion}.
For sufficiency, we give an algorithm to construct a URVG layout:

{\bf INPUT:} A partition of a tree $T$ into two caterpillar forests
$F_1$ and $F_2$, each having maximum degree 3.

{\bf OUTPUT:} A URVG layout of $T$ with $T_X = F_1$ and $T_Y=F_2$.

Choose an arbitrary root $r$ for $T$ and give $T$ a breadth-first numbering,
$v_1, \ldots, v_n$, starting with $r=v_1$. For each $i$, let $S_i$ denote
the square in the URVG layout representing $v_i$, and let $(x_i, y_i)$ be the
coordinates of the lower left corner of  $S_i$. The algorithm places square $S_1$
arbitrarily, and then, for $i = 1, \ldots, n-1$, it places squares representing
the children of vertex $v_i$. For each $i = 1, \ldots, n-1$, at the point in the
algorithm at which we have placed a square for $v_i$ and squares for
all its children, but {\em no} higher-numbered squares, the following set of
invariants is maintained.
\smallskip

{\bf Algorithm Invariants:}
\begin{enumerate}%i1
  \item\label{hflush} A leg edge in $F_1$ corresponds to a horizontal flush
  visibility in the layout.
  \item\label{vflush} A leg edge in $F_2$ corresponds to a vertical flush
  visibility in the layout.
  \item\label{hprot} A spine edge in $F_1$ corresponds to a horizontal
  protruding (i.e., not flush) visibility in the layout.
  \item\label{vprot} A spine edge in $F_2$ corresponds to a vertical
  protruding visibility in the layout.
  \item\label{onlymom} If a vertex $v_j$ has parent $v_i$ in the
  breadth-first numbering, where $i < j$, then at the point when the
  square $S_j$ is placed, it sees no square other than $S_i$.
\end{enumerate}%i1

We observe that if $S_i$ is a square in the layout corresponding to vertex $v_i$, then $v_i$ may have either at most  two incident leg edges in $F_1$, or at most one incident leg edge and two incident spine edges in $F_1$. In the former case the algorithm places the squares for the two leg neighbors of $S_i$ flush with it on opposite sides. In the latter case, the square for the leg neighbor is placed flush on one side of $S_i$, and both squares for the spine neighbors are placed on the opposite side of $S_i$, protruding so that they both see $S_i$ but not one another. A similar statement holds for the neighbors of $v_i$ in $F_2$; see Fig.~\ref{fig:treeneighbors}. The algorithm chooses the exact placement of each square so as not to introduce unwanted visibilities with other squares that have already been placed.

\begin{figure}[htb]%[htbp]
 \centerline{\epsfig{file=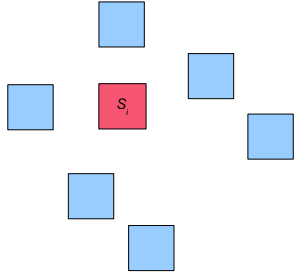, width=0.3\textwidth}}
 \caption{Placement of neighbor squares in a tree layout.}
 \label{fig:treeneighbors}
\end{figure}

{\bf Algorithm to create a URVG layout of the tree $T$:}
\begin{enumerate}%e1
  \item Let $i = 1$, and arbitrarily place square $S_1$ representing $v_1$.  Define $x_m$ to be the smallest $x$-coordinate of any square that has been placed in the layout and $x_M$ to be the largest $x$-coordinate of any square that has been placed in the layout.  Similarly, define $y_m$ and $y_M$ to be the smallest and largest $y$-coordinates of any squares that have been placed in the layout.
  \item Place squares representing the children of $v_1$ as follows, maintaining the invariants:
    \begin{itemize}%i2
    \item There are at most two leg edges of $F_1$ incident with $v_1$. Place squares representing the endpoints of these edges, if they exist, at position $(x_1-2, y_1)$ for the first (in breadth-first order) and $(x_1+2, y_1)$ for the second. For leg edges of $F_2$ incident with $v_1$, place the squares at $(x_1, y_1-2)$ and $(x_1, y_1+2)$.
    \item If $v_1$ is incident to any spine edges in $F_1$, then it is incident to at most two such edges, and it is incident to at most one leg edge of $F_1$, already placed at $(x_1-2, y_1)$. Place squares representing the endpoints of the spine edges, if they exist, at $(x_1+2, y_1+2/3)$ and $(x_1+4, y_1-2/3)$, in breadth-first order. For spine edges of $F_2$ incident with $v_1$, place the squares at $(x_1+2/3, y_1+2)$ and $(x_1-2/3, y_1+4)$.
    \end{itemize}%12
  \item\label{item_inc} Let $i = i+1$. Square $S_i$ has been placed as a child
  of its parent square $S_p$, where $1\le p<i$. If $i = n$, the layout is
  complete, since $v_n$ has no children. Otherwise, we place squares for
  each child $v_j$ of $v_i$ as follows:
    \begin{itemize}%i3
    \item Leg edges to children of $v_i$: If $v_i$ has a leg edge $v_i$--$v_j$ in $F_1$, incident to a child $v_j$,
    the placement of $S_j$ depends on the relative placements of $S_i$ and its
    parent square $S_p$.  If $v_i$--$v_p$ is also an edge (leg or spine) of $F_1$, assume without loss of generality (since we can flip the existing layout horizontally or vertically if necessary) that its square was placed to the {\em left} of $S_i$. Then we place the square $S_j$ to the {\em right} of $S_i$, at position $(x_M+2, y_i)$. Otherwise we place the (at most two) leg edges of $F_i$ from $v_i$ to its children, in breadth-first order, at positions $(x_m-2, y_i)$ and $(x_M+2, y_i)$. We place squares representing the endpoints of (at most two) leg edges of $F_2$ from $v_i$ to its children in an analogous manner, using positions $(x_i, y_m-2)$ and $(x_i, y_M+2)$. The values of $x_m$, $x_M$, $y_m$, and $y_M$ are updated after the placement of each square.
        \item Spine edges to children of $v_i$: If $v_i$ has a spine edge $v_i$--$v_j$ in $F_1$, incident to a child $v_j$,
    the placement of $S_j$ depends on the relative placements of $S_i$ and its parent square $S_p$.  If $v_p$--$v_i$ is a leg edge of $F_1$, then $S_p$ is a horizontal, flush neighbor of $S_i$; assume without loss of generality that $S_p$ is flush neighbor lying to the {\em right} of $S_i$. If $v_p$--$v_i$ is a spine edge of $F_1$, then $S_p$ is a horizontal, protruding neighbor of $S_i$. Without loss of generality, we assume that $S_p$ lies to the {\em left} of $S_i$ and is a downward protruding neighbor of $S_i$. 

We wish to place
        $S_j$ as an upward protruding neighbor of $S_i$, also lying to the left of
        $S_i$ (on the side opposite of $S_p$ if $v_p$ is a leg neighbor of $v_i$ in $F_1$, and on the same side if $v_p$ is a spine neighbor of $v_i$ in $F_1$). Invariant \ref{onlymom} guarantees
        that a horizontal line through the upper edge of $S_i$ intersects the interiors
        of no other squares in the layout. We replace that line with a horizontal band
        of height one to create room for an upward protruding neighbor of $S_i$; see
        Fig.~\ref{fig:makeroom}. If $x_m$ is the smallest $x$-coordinate in the current
        layout, then the coordinates of $S_j$ are $(x_m-2, y_i+2/3)$, and the
        invariants are maintained. The values of $x_m$, $x_M$, $y_m$, and $y_M$ are updated after placement. \smallskip 

If $v_p$--$v_i$ is not a spine edge of $F_1$, then place the square for the first spine neighbor in $F_1$ to the left of $S_i$ as described above. To place a square for a second spine neighbor of $v_i$ in $F_1$, we replace the line through the lower edge of $S_i$ with a horizontal band of height one, creating room to place the square at position $(x_m-2, y_i-2/3)$. The values of $x_m, x_M, y_m$, and $y_M$ are updated after placement. \smallskip

We place squares representing the endpoints of (at most two) spine edges of $F_2$ from $v_i$ to its children in an analogous manner, using positions $(x_i-2/3, y_m-2)$ and $(x_i+2/3, y_m-2)$. 
%The values of $x_m$, $x_M$, $y_m$, and $y_M$ are updated after the placement of each square.
        \end{itemize}%i3
    \item At this point squares representing all the children of $v_i$ have
    been placed. Return to step \ref{item_inc}.
\end{enumerate}
\end{proof}

\begin{figure}[htb]%[htbp]
 \centerline{\epsfig{file=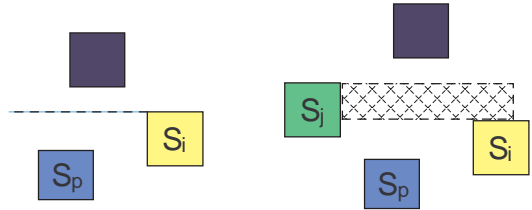, width=0.6\textwidth}}
\caption{Making room for a child square $S_j$ of $S_i$.}
\label{fig:makeroom}
\end{figure}

\begin{figure}[htb]
  \centerline{\epsfig{file=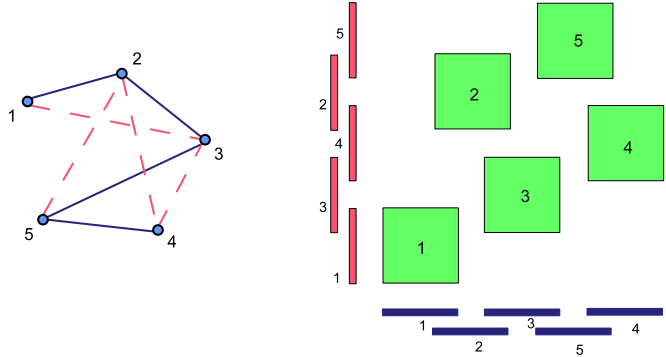, width=0.7\textwidth}}
  \caption{URV layout of a graph with linear arboricity 2.}
\label{fig:linarb2}
\end{figure}

A comparable result holds for general graphs, if we require a decomposition
into two forests of \emph{paths}, rather than into subdivided caterpillar forests
with maximum degree 3. The {\em linear arboricity} of a graph $G$ is the minimum
number of linear forests whose union is $G$. Similarly the {\em caterpillar arboricity}
of $G$ is the minimum number of caterpillar forests whose union is $G$. It is shown in
(\cite{BDHS97}, Theorem 5)
that, if $G$ has caterpillar arboricity 2, then it is an RVG. The proof
is constructive: each caterpillar forest is represented as an interval graph, one
along the $x$-axis and the other along the $y$-axis. The Cartesian product of
horizontal and vertical intervals corresponding to the same vertex is a rectangle in the plane, and the resulting set of
rectangles is an RVG representation of $G$. If it happens that both caterpillar
forests are actually linear forests, then  the intervals can all have equal length,
making $G$ a URVG; an example is shown in Fig.~\ref{fig:linarb2}. Hence the
following result follows immediately.

\begin{thm}\label{thm:linarb2}
If $G$ has linear arboricity 2, then $G$ is a URVG.
\end{thm}

\begin{figure}[htb]
  \centerline{\epsfig{file=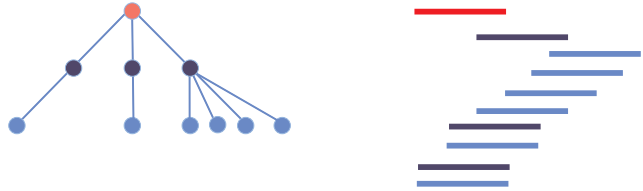, width=0.7\textwidth}}
  \caption{Weak UBV layout of a tree.}
\label{fig:weakTree}
\end{figure}

The next result states that, in contrast to Theorem \ref{thm:URVGtree}, 
every tree is a weak UBVG.  See Fig.~\ref{fig:weakTree} for a sample layout.

\begin{thm}\label{thm:weakURVGtree}
Every tree is a weak UBVG, and hence a weak URVG.
\end{thm}

\begin{proof}
We prove this by providing an algorithm to construct a weak 
UBV layout.
Let $T$ be a tree.  Choose an arbitrary root $r$, and give $T$ a 
breadth-first numbering starting with $r=0$.
We specify the position of a bar corresponding to vertex $v$ by the 
coordinates of its left endpoint $(v_x, v_y)$.
We begin by placing the bar corresponding to $r$ at position $(0,0)$.
If $r$ has $k_0$ children, place the bar for child $i$ of $r$ at 
position $(1 - i/k_0, -i)$.
Next we will place the vertices on level two.
If parent $p$ has $k_p$ children insert a horizontal band of height
$k_p$ below the bar corresponding to $p$ in the layout.
Child $i$ of parent $p$ should then be placed at position $(x_p+1-i/k_p, y_p-i)$.
Continue in this manner to place the children on the remaining levels.
\end{proof}

\section{Complete Bipartite Graphs}\label{sec:bipartite}

In this section we establish which complete bipartite graphs are URVGs,
and which are weak URVGs.
Throughout, we always write $K_{m,n}$ with $m\le n$,
and we denote by $V_m$ and $V_n$ the two partite sets of $K_{m,n}$.

\begin{thm}\label{thm:kmnsuff}
The complete bipartite graph $K_{m,n}$, $m\le n$, is a URVG if $m\le2$
and $n\le 6$, or $m=3$ and $n \leq 4$. $K_{m,n}$ is a weak URVG if $m\le2$ (and $n$
is arbitrary), or $m\le3$ and $n\le4$.
\end{thm}

\begin{figure}[htb]%[htbp]
  \centerline{\epsfig{file=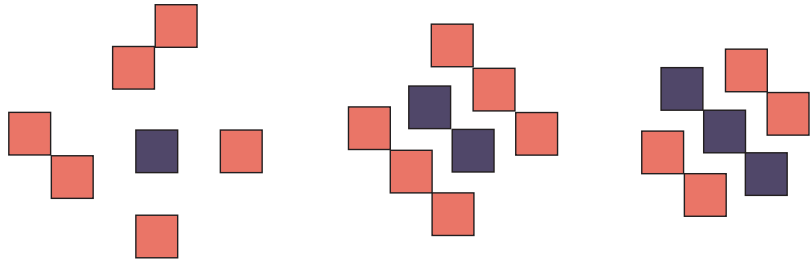, width=0.9\textwidth}}
\caption{URV Layouts of $K_{1,6}$, $K_{2,6}$, and $K_{3,4}$. }
\label{fig:kmnlayouts}
\end{figure}

\begin{figure}[htb]%[htbp]
  \centerline{\epsfig{file=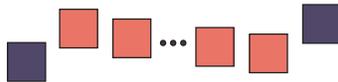,  width=0.4\textwidth}}
\caption{A weak URV layout of $K_{2,n}$.}
\label{fig:k2nweak}
\end{figure}

\begin{proof}
URV layouts of $K_{1,6}$, $K_{2,6}$, and $K_{3,4}$ are shown in Fig.~\ref{fig:kmnlayouts}.
A weak URV layout of $K_{2,n}$ is shown in Fig.~\ref{fig:k2nweak}. Layouts for smaller
values of $n$ are obtained by deleting squares from these layouts.
\end{proof}

The results of the remainder of this section prove that the conditions
of Theorem~\ref{thm:kmnsuff} are also necessary. 
%We also show that $K_{2,n}$ and $K_{3,4}$
%are the largest complete bipartite graphs having weak URV layouts.

\begin{figure}[htb]%[htbp]
  \centerline{\epsfig{file=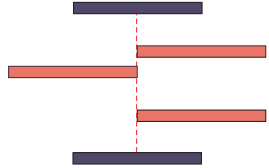,  width=0.3\textwidth}}
\caption{Two-sided UBV layout of an empty cycle.}
\label{fig:emptycycle}
\end{figure}

Given a URV layout of $K_{m,n}$, we write $K_X$ and $K_Y$ 
for $(K_{m,n})_X$ and $(K_{m,n})_Y$.  We call
a cycle in a plane graph \emph{empty} if it bounds a finite face. As noted in 
(\cite{DV03}, Proposition 1.10),  
the layout of an empty $n$-cycle $C$ in a UBVG corresponds to a vertical line $\ell_C$
joining the interiors of the top and bottom cycle bars, such that each of the other $n-2$
`intermediate' bars of $C$ has its left or right endpoint on $\ell_C$. See
Fig.~\ref{fig:emptycycle}.
We call a UBV layout of an empty cycle \emph{one-sided} if all the intermediate
bars lie on the same side of $\ell_C$ and \emph{two-sided} otherwise.

\begin{lem}\label{lem:emptycycle}
If $C$ is a cycle in a UBVG, then in the plane embedding induced by
the corresponding UBV layout, every edge of $C$ lies on an empty cycle.
\end{lem}

\begin{figure}[htb]%[htbp]
  \centerline{\epsfig{file=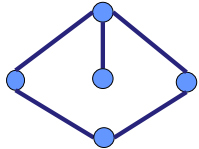,  width=0.2\textwidth}}
\caption{A cycle whose edges do not lie on an empty cycle.}
\label{fig:nonemptyface}
\end{figure}

\begin{proof}
It follows from the characterization of BVGs in (\cite{Wismath85},Theorem 4) and (\cite{TT86}, Theorem 5) that the plane embedding induced by a UBV layout has all cutpoints on the exterior face, but we must show that no 2-connected block incident with such a cutpoint lies in the interior of a finite face; see Fig.~\ref{fig:nonemptyface}.
We claim that the subgraph $G_C$ comprised of the vertices
and edges lying on $C$ and in its interior is 2-connected. First,
$C$ itself is a 2-connected graph. Next, let  $b_w$ be a bar in the interior
of $G_C$. The set comprised of the bars of $C$, together with vertical lines
representing visibilities corresponding to the edges of $C$, contains a Jordan
curve. Thus, if we pass a vertical line through $b_w$, it will induce a path of
visibilities containing $b_w$  from one vertex $u$ of $C$ to another vertex $v$
of $C$.  Hence, every vertex in $G_C$ lies on a path joining two vertices of  $C$.
Therefore $G_C$ is 2-connected, and so its internal faces are bounded by empty cycles.
\end{proof}

\begin{lem}\label{lem:cyclelemma}
If there is a URV layout of $K_{m,n}$ with an empty cycle $C$ in $K_Y$,
then $C$ is a 4-cycle, $m=2$, and the maximum degree 
of any vertex in $K_Y$ is 3. Furthermore,
\begin{enumerate}
  \item If $C$ is one-sided, then $2 \leq n \leq 3$;
  \item If $C$ is two-sided, then $2 \leq n \leq 4$.
\end{enumerate} An analogous result holds when $C$ is in $K_X$.
\end{lem}

\begin{figure}[htb]%[htbp]
  \centerline{\epsfig{file=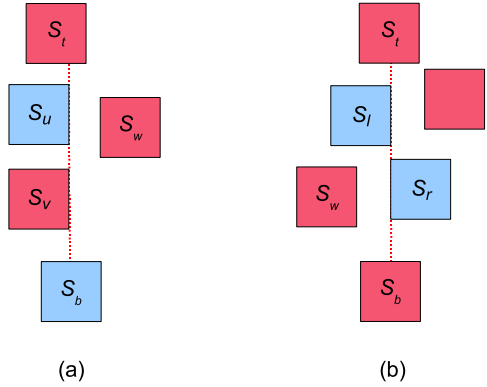,  width=0.5\textwidth}}
\caption{The two cases of Lemma~\ref{lem:cyclelemma}.}
\label{fig:cyclelemma}
\end{figure}

\begin{proof}
First we note that $C$ has at most two left-intermediate and 
two right-intermediate bars: if $C$ has, say, three or more left-intermediate 
bars, let $Y_p$ be the third-highest one of these. Then $p$ and $t$ 
are in different partite sets, but $S_p$ cannot see $S_t$. Next note that there 
cannot be both two left-intermediate bars {\em and} two right-intermediate bars, 
because then $t$ and $b$ would be in different partite sets, but $S_t$ could 
not see $S_b$. So $C$ is a 4-cycle, and either $C$ has one left-intermediate 
and one right-intermediate bar or, without loss of generality, two 
left-intermediate bars and no right-intermediate bars.

Suppose next that $C$ has two left-intermediate bars with $Y_u$ the higher of the two. 
Without loss of generality, assume that $x_t \le x_b$, as illustrated in 
Fig.~\ref{fig:cyclelemma}(a). Then the only way there can be additional squares 
is if $x_t <x_b$, and a square $S_p$ is placed with $x_t+1\le x_p < x_b$, so 
that $S_p$ sees $S_u$ horizontally and $S_b$ vertically. Furthermore at most one 
such square can be placed this way, thus $m=2$ and $2\le n\le 3$.

Lastly, assume  $C$ has one left-intermediate and one right-intermediate bar. 
Then $t$ and $b$ are in the same partite set, and the intermediate bars prevent 
them from having any other common neighbors in the URV layout. Hence $m=2$. 
It is possible to place at most two more squares in the same partite set as $t$ 
and $b$, as illustrated in Fig.~\ref{fig:cyclelemma}(b), so that $2 \le n \le 4$. 
Furthermore, each of these squares must see one intermediate square horizontally 
and the other vertically, so that the maximum degree 
of any vertex in $K_Y$ is 3.
\end{proof}

\begin{lem}
If $n\ge7$, then $K_{m,n}$ is not a URVG.
\end{lem}

\begin{proof}
We assume $K_{m,n}$ has a URV layout and obtain a contradiction. Since $n\ge7$,
there is a vertex $s$ in, say, $K_V$, that has vertical degree $\ge4$. By 
(\cite{DV03}, Corollary 1.7)
the vertex $s$ lies on a cycle in $K_V$, and so by Lemma~\ref{lem:emptycycle}, $s$
lies on an empty cycle. This contradicts the degree conditions of Lemma~\ref{lem:cyclelemma}.
\end{proof}

\begin{thm}\label{thm:weakKmn}
If $m,n\ge4$, or $m\ge3$ and $n\ge 5$, then $K_{m,n}$ is not a weak URVG.
\end{thm}

\begin{proof}
Let $m,n\ge4$ or $m\ge3, n\ge5$, and suppose that $L$ is a weak URV layout of
$K_{m,n}$.
%Given a unit square $S_u\in L$, remember that $x_u$ and $y_u$ denote the
%coordinates of its lower left corner.
%, and for convenience, let $x'_U=x_U+1$,
%$y'_U=y_U+1$ be the coordinates of its upper right corner.
If necessary, we modify the layout slightly so that it becomes noncollinear
without losing any visibilities. By a {\em monotonic set of squares} we mean
a set of squares whose set of $y$ coordinates forms a monotonic sequence
when listed in increasing order of the $x$ coordinates. The conditions on $m,n$
imply that there are eight or more vertices, so by Theorem~\ref{thm:subseq}  there is a monotonic set of three squares in the layout $L$. We assume without
loss of generality that this set is monotonically increasing. We consider two cases
depending on whether or not such a sequence exists with all elements in the same
partite set.

{\bf Case 1:} There is an increasing set of three squares
$\{S_1,S_2,S_3\}\subseteq V' = V_m$ or $V_n$. Let $V''=V(K_{m,n})\setminus V'$, and
assume without loss of generality the squares are 
labeled from left to right, $S_1,S_2,S_3$. Note that no square $S_b$ with 
$b\in V''$ can be part of a larger increasing sequence
containing $S_1, S_2, S_3$, because then $S_b$ would be able to see at most two other squares
in the sequence. For an element $b\in V''$, consider the 3-tuple of directions
from which $S_b$ sees $S_i, i = 1, 2, 3$. For example, the 3-tuple $(N,N,W)$
signifies that $S_b$ is above $S_1$ and $S_2$, so sees them from the north, 
and is to the left of $S_3$, so sees it from the west.
Because $S_1, S_2, S_3$ is an increasing sequence, a 3-tuple must have at least
two consecutive repeated terms. The triples $(N,N,E)$, $(W, S,S)$,
$(E,E,N)$, and $(S,W,W)$ are prohibited.
The remaining possible 3-tuples are:
$$(N,N,N), (N,N,W),(N,N,S),  (S,S,S), (E,S,S), (N,S,S),  $$
$$(E,E,E), (E,E,S), (E,E,W),(W,W,W), (N,W,W), (E,W,W).$$

Since $S_1, S_2, S_3$ is increasing, if $S_b$ sees two squares from the 
north or east they must be $S_1$ and $S_2$, and if from the south or west
they must be $S_2$ and $S_3$.  Note that when $S_b$ sees  two squares
from the same direction, it blocks any other square from simultaneously seeing 
those two squares from that direction, so $|V''| \leq 4$.

We next show that, while $V''$ may have three or four elements, 
$V'$ can have only three, so it is not the case that $m,n \geq 4$ or $m \geq 3$
and $n \geq 5$.

We label the four possible elements that $V''$ could have by $S_n,S_s,S_e,$
$S_w$, to indicate the direction from which they each see more than one of
$S_1,S_2,S_3$, as prescribed by the 3-tuples above. The following argument applies
whether or not $S_w$ is present (indicated by [\ ]), 
and by symmetry it applies if any three of the
four are present. Consider the square $S_2$, which all four of these squares see in the
directions prescribed by their subscripts. There are four disjoint regions of the
plane composed of points not visible to $S_2$, located to the northwest, southwest,
southeast, and northeast of $S_2$; see Fig.~\ref{fig:corridors}(a). Because $S_n$ sees
both $S_1$ and $S_2$ from the north, and because $\{S_1,S_2,S_3\}$ is increasing,
$S_n$ must intersect the region northwest of $S_2$. Similarly, [$S_w$ also intersects
this region, while] $S_s$ and $S_e$ both intersect the region southeast of $S_2$.
If there is an additional square $S_4\in V'$ in the layout, it cannot intersect the
region northwest of $S_2$, since then it cannot see either $S_s$ or $S_e$; likewise
it cannot intersect the region southeast of $S_2$ since it must see $S_n$ [and $S_w$].
Assume without loss of generality that $S_4$ intersects the region northeast of $S_2$
(recall our assumption that the layout is noncollinear). Since both $S_1,S_2,S_3$ and
$S_1, S_2,S_4$ are increasing sets, we relabel $S_3$ and $S_4$, if necessary, so that
$S_3$ is further left than $S_4$. Now $S_s$ sees $S_3$ and  $S_4$ from the south,
so $S_e$ cannot see both these squares from the south. Hence $S_e$ must see the one
that is further left, namely $S_3$, from the east, forcing $S_3$ to intersect the
visibility corridor to the east of $S_2$. But this prevents $S_4$, which intersects
the region northeast of $S_2$, and which is further right than $S_3$, from seeing $S_s$.
This conclusion is reached whether or not $S_w$ is present, so by symmetry, we conclude
that $V'$ has at most three elements.

\begin{figure}[htb]
  \centerline{\epsfig{file=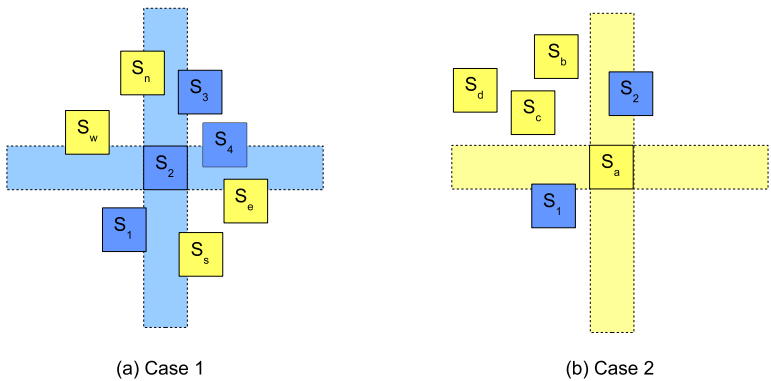, width=0.8\textwidth}}
\caption{Positions of squares in proof of Theorem~\ref{thm:weakKmn}.}
\label{fig:corridors}
\end{figure}

{\bf Case 2:} There is no monotonic subsequence of length 3 or more with all
elements in one of $V',V''$. By Theorem~\ref{thm:subseq}, this can occur only
if both $m$ and $n$ are $\le 4$ so $m=n=4$. 
Again by Theorem~\ref{thm:subseq},  there is a
strictly increasing subsequence of length 3 or more with two elements in $V'$ and
one element in $V''$. Note  that the element of $V''$ in the increasing sequence of
three must be the middle  in order to see both of the other two, so we name these
elements left to right, $S_1, S_a, S_2$. As in Case 1, consider the regions
northwest, southwest, southeast, and northeast of the middle square, $S_a$. The square
$S_1$ can intersect only the southwest region, and $S_2$ can intersect only the
northeast region. No other element of $V''$ can intersect these two regions, because then it
cannot see whichever of $S_1$ or $S_2$ intersects only the diagonally opposite region.  Without
loss of generality, assume that there is another square $S_b$ 
that intersects the region northwest of $S_a$.
Then no element  of $V''$ can intersect the southeast region, because  it, together
with $S_a$ and $S_b$, would form a monotonic set of length 3, contradicting the assumption
of Case 2.  Suppose there are three or more elements of $V''$ in the region northwest of
$S_a$; see Fig.~\ref{fig:corridors}. 
Since none can have y-coordinates less than $S_1$'s and they can't be monotonic, then 
some two out of the three form a monotonic set with $S_1$, giving a contradiction.
Thus $V''$ has at most three elements.
\end{proof}

\section{Edge Bounds}\label{sec:bounds}       % Enter section title between curly braces

In this section we give several results bounding the number of edges of  URVGs. First we
use the characterizations of Theorems~\ref{thm:UBVGtree} and \ref{thm:URVGtree} to give
tight upper bounds on the number of edges in a depth-$s$ UBV tree and in a depth-$s$ URV tree. Next
we modify the methods used for RVGs in \cite{HSV99} and \cite{DH94} to bound the number of
edges in general URVGs and bipartite URVGs, respectively, and we give examples to show that these
bounds have tight order. \\

%\subsection{%Edge bounds for depth-$s$ trees}\label{sub:treebounds}
\noindent \textbf{Edge bounds for depth-$s$ trees} 

\begin{thm} \label{thm:UBVtreebound}
If $T$ is a  rooted depth-$s$ tree that is a unit bar visibility graph, then $T$  has at
most $s^2+2s$ edges.
\end{thm}

\begin{figure}[htb]
  \centerline{\epsfig{file=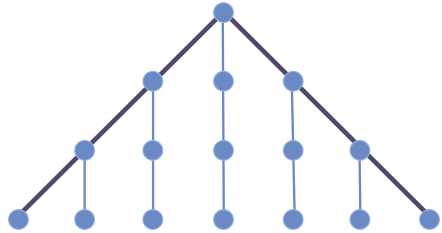, width=0.5\textwidth}}
\caption{The depth-3 UBVG tree $T_{B,3}$.}
\label{fig:UBVGTree}
\end{figure}

\begin{proof}
Throughout we use the result of Theorem~\ref{thm:UBVGtree}, that a tree is a UBVG if
and only if it is a subdivided caterpillar forest with maximum degree 3. Define
$T_{B,s}$ to be the subdivided caterpillar whose spine has length $2s+1$, rooted at the
center spine vertex, and with a leg at each vertex extending to depth $s$, as illustrated in
Fig.~\ref{fig:UBVGTree}. By Theorem~\ref{thm:UBVGtree}, $T_{B,s}$ is a UBVG, and it's
easy to see that the number of edges of $T_{B,s}$ equals $3 + 5 + \ldots + (2s+1) =
\sum_{k=1}^{s}(2k+1)=s^2+2s$. If $T\ne T_{B,s}$ and $T$ is rooted at a spine vertex,
then $T$ is a subtree of $T_{B,s}$. If $T\ne T_{B,s}$ and $T$ is \emph{not} rooted at
a spine vertex, suppose that $v$ is the spine vertex closest to the root $r$ of $T$,
and that $v$ is a depth-$q$ vertex, where $1 < q \le s$. Then the subtree rooted at $v$
is a subtree of $T_{B,s-q}$ with no leg at $v$, so it has at most $(s-q)^2+2(s-q)-(s-q)$
$=s^2-2sq+q^2+s-q$ 
edges. The path from $v$ to $r$ adds another $q$ edges plus $s$ possible edges 
on a path from $r$ to depth $s$, for a total of at most
$s^2 + 2s - 2 s  q + q^2 =(s-q)^2+2s <s^2+2s$ edges.
\end{proof}

\begin{thm} \label{thm:URVtreebound}
Let $T$ be a rooted URV tree.  Then the number of vertices at depth $k \geq 1$
is bounded by $c_k$, where $c_k$ is given recursively by the linear recursion relations
$c_k=a_{k-1}+4c_{k-1}$ and $a_k=a_{k-1}+2c_{k-1}$ with initial values $a_1=4$ and $c_1=6$.
Furthermore, this bound is tight. 
%If $T$ is a  rooted depth-$s$  tree that is a unit rectangle visibility graph,
%then the maximum number of vertices at depth-$s$ in $T$ is given by the
%coefficient of $x^s$ in the Taylor series for the function
%$V(x) = \frac{2x(3-x)}{2x^2-5x+1}+1$.
\end{thm}

\begin{proof}
%Throughout we use the result of Theorem~\ref{thm:URVGtree}, that a tree $T$
%is a URVG if and only if it is the union of two subdivided caterpillar forests,
%each with maximum degree 3. The result is proved in two parts. 
%First, in a manner analogous to the %definition of $T_{B,s}$
%in the preceding proof, we define a depth-$s$ tree $T_{R,s}$ that is a union of
%two subdivided caterpillar forests, and we show that it satisfies the condition of the theorem. 
We first exhibit a canonical depth-$s$ URV tree $T_{R,s}$ achieving this bound,
as this motivates the given recursion relations.  The tree $T_{R,s}$ is defined
analogously to $T_{B,s}$ in the proof of Theorem~\ref{thm:UBVtreebound}.
Its root $r$ is the center vertex of two spines, each of length
$2s+1$, for two subdivided, maximum degree 3 caterpillars, one in each forest.
We distinguish the two forests by calling one \emph{red} and 
the other \emph{blue}. For each vertex $v$
from level 1 to level $s-1$, the number of $v$'s children and the colors (red or blue)
and types (spine or leg) of  the edges from $v$ to its children are determined by the
edge $e$ from $v$ to its parent. If $e$ is a blue spine edge, then $v$ has five children,
one incident with a blue spine edge, two with red spine edges, one with a blue leg edge,
and one with a red leg edge; a symmetric definition holds if $e$ is red spine edge.
If $e$ is blue leg edge, then $v$ has four children, two incident with red spine edges,
one with a blue leg edge, and one with a red leg edge; a symmetric definition holds
if $e$ is red leg edge. In this case, the current subdivided caterpillar containing
the edge $e$ is continued at the vertex $v$, and a new subdivided caterpillar of
the opposite color begins with $v$ as its root. Part of the tree $T_{R,3}$ is
shown in Fig.~\ref{fig:URVGTree}  (spine edges are thicker lines, and leg edges
are thinner lines).

\begin{figure}[htb]
  \centerline{\epsfig{file=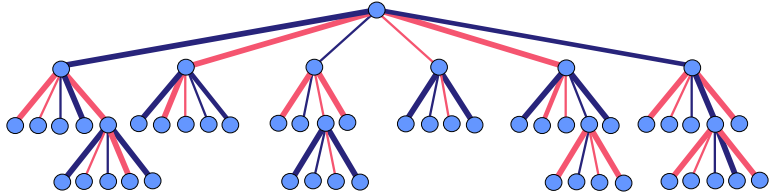, width=0.8\textwidth}}
\caption{Part of the depth-3 URVG tree $T_{R,3}$.}
\label{fig:URVGTree}
\end{figure}

We now count the vertices on each level $k$ of $T_{R,s}$. For $1\le k\le s-1$,
let $a_k$ be the number of vertices at level $k$ with five children, i.e., those with spine edges to their parents; let $b_k$
be the number of vertices at level $k$ with four children, i.e. those with leg edges to their parents; and let $c_k=a_k+b_k$
be the total number of vertices on level $k$. 

Note that  each vertex $v$ at level $k$ with five children has three
children who themselves have five children; the other two children each have  four
children. 
If $v$ has four children, then two of them have five children and two of
them have four children. We therefore have the following recurrence equations:
\begin{eqnarray}
{
a_k=3 a_{k-1} + 2 b_{k-1}, 
b_k=2 a_{k-1} + 2 b_{k-1},
a_1=4, b_1=2.}
\label{eqn:akbk}
\end{eqnarray}

Using the fact that $c_k=a_k+b_k$, we can eliminate the $b_k$ term to get a system
involving $a_k$ and $c_k$:
\begin{eqnarray}
{
a_k=a_{k-1} + 2 c_{k-1}, %\\
c_k=a_{k-1} + 4 c_{k-1},
a_1=4, c_1=6.}
\label{eqn:akck}
\end{eqnarray}

\medskip
Now let $T$ be any depth-$s$ URV tree, with a decomposition into 
two UBV forests (one red and one blue) as guaranteed by Theorem~\ref{thm:URVGtree}.
We say $T$ is a \emph{downward} tree if every leg extends 
strictly downward from its point of attachment to its spine. If $T$ is a downward tree, then $T$ is a subtree 
of $T_{R,s}$, and in fact can be mapped onto $T_{R,s}$ preserving the UBV decomposition.
This can be seen easily by induction on $s$.  Thus for every downward tree, 
the bound holds.

Now suppose $T$ has an \emph{up-leg}, that is, a leg that does not extend strictly
downward from its spine.  We observe the following consequences of $T$ being a rooted tree: 

\begin{itemize}
  \item Each caterpillar in the decomposition may have at most one such leg.
  \item This leg must be attached to the highest vertex of the spine.
  \item This leg must extend strictly upward to its highest vertex and then strictly downward from there.
\end{itemize}

Consider each vertex $x$ of $T$ in breadth first order.
If an up-leg intersects its spine at vertex $x$, we perform 
the following surgery on $T$.  
We convert the portion of the spine that is to the left of $x$ to a leg, 
and convert the up-leg into a spine, extending its terminus downward if necessary 
until it reaches at least the same depth as the old left side of the spine.
We then remove the edge connecting each leg to the old spine, and 
add an edge connecting the leg to the vertex of the new portion of the 
spine at precisely the same depth as the one it was attached to on the old spine.
All other edges of $T$ remain unchanged.  See Figures~\ref{fig:uconfigs} 
and \ref{fig:ureplace}.  

\begin{figure}[htb]%[htbp]
  \centerline{\epsfig{file=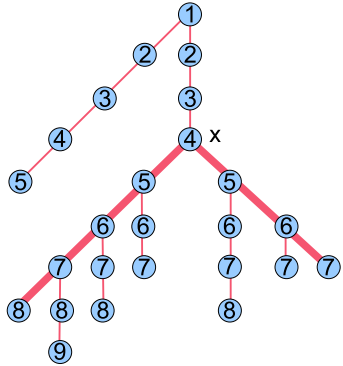, width=0.35\textwidth}}
  \caption{A tree with an up-leg adjacent to vertex $x$.  The number on each vertex 
corresponds to its depth in the tree.}
  \label{fig:uconfigs}
\end{figure}

\begin{figure}[htb]%[htbp]
  \centerline{\epsfig{file=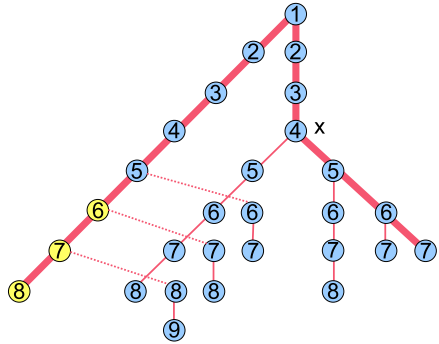, width=0.45\textwidth}}
  \caption{The modified tree, removing the up-leg from Fig.~\ref{fig:uconfigs}.
The new edges are represented by dashed lines.}
  \label{fig:ureplace}
\end{figure}

The result is a downward URV tree with at least as many vertices at each depth 
as the original tree, and hence no more than $c_k$ vertices at level $k$.

\end{proof}

Standard linear recursion techniques can be used to compute 
$c_k$, either  as the coefficient of $x^k$ in the generating function 
$V(x) = \sum_{k=0}^{\infty}{c_kx^k}$ $=\frac{2(3-x)}{2x^2-5x+1}+1$, 
%The first six values of $c_k$ are $c_0=1, c_1=6, c_2=28, c_3=128, c_4=584, c_5=2,664$. 
or in the following closed form:
\begin{align}
 c_k=&\frac{1}{17} 2^{-1 - 2 k} ((17 - 7 \sqrt{17}) (10 - 2 \sqrt{17})^k +\label{eqn:cs}\\
  &2^{1 + k} \sqrt{17} (5 + \sqrt{17})^k +  (2 (5 + \sqrt{17}))^
    k (17 + 5 \sqrt{17})).\notag
\end{align}

\begin{cor}
If $T$ is a  rooted depth-$s$ URV tree,
then the maximum number of edges in $T$ is
\begin{align}
-2+&\left(1-\frac{3}{\sqrt{17}}\right) \left(\frac{1}{2} \left(5-\sqrt{17}\right)\right)^s+\label{eqn:edges}\\
&\left(1+\frac{3}{\sqrt{17}}\right) \left(\frac{1}{2} \left(5+\sqrt{17}\right)\right)^s\notag
\end{align}
\end{cor}

\begin{proof}
This is the sum of Equation (\ref{eqn:cs}) as $k$ ranges from 1 to $s$.
\end{proof}
 
%Alternatively, we can compute the sum of Equation (\ref{eqn:cs}) from 1 to $s$ and simplify, to obtain %the following expression for the maximum number of  edges in a depth-$s$ URV tree:

%\subsection{Edge bounds for general URVGs}\label{sub:generalbounds}
\noindent \textbf{Edge bounds for general URVGs} 

\begin{thm}\label{thm:bound}
For $n\ge 1$, let $G$ be a URVG with $n$ vertices. Then $|E(G)| \le 6n-4\lceil\sqrt{n}\rceil +1$.
\end{thm}

\begin{figure}
  \centerline{\epsfig{file=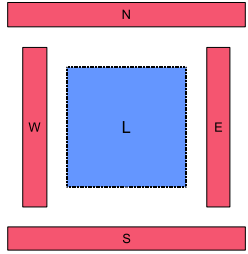, width=.25\textwidth}}
\caption{A URV layout plus 4 surrounding rectangles.}
\label{fig:enclosedlayout}
\end{figure}

\begin{proof}
For $n=1, \ldots, 10$ this bound follows immediately from the edge sets
of complete graphs. 
Let $L$ be a URV layout of $G$. Surround the layout $L$ with four rectangles (that are
not unit squares) labeled $N, S, E$, and $W$, as shown in Fig.~\ref{fig:enclosedlayout}.
Call the induced rectangle visibility graph $G_+$. Partition the edges of $G_+$ into the
two sets $E = E(G)$ and $E' = E(G_+)-E$. Hence $E'$ comprises the edges of $G_+$ that
have either one or two endpoints in the set $\{N, S, E, W\}$. By a result of (\cite{HSV99},
Theorem 1),
any rectangle visibility graph with $p \ge 5$ vertices has at most $6p-20$ edges.
Therefore,
$$
|E| + |E'| = |E(G_+)| \le 6(n+4) - 20 = 6n+4
$$

Now note that the number of edges of $E'$ with exactly one end point in 
$\{N,S,E,W\}$ is at least as large as the number of squares on the perimeter of a rectangle
containing all $n$ squares of G.  This is at least $\lceil4\sqrt{n}\rceil$, which is greater than 
or equal to  $4\lceil\sqrt{n}\rceil -3$.  There are 6 edges with both endpoints 
in $\{N,S,E,W\}$, so 
$$
|E'| \ge 
(4\lceil\sqrt{n}\rceil -3)+6 = 
4\lceil\sqrt{n}\rceil +3.
$$
Therefore,
$$
|E| \le
6n + 4 - |E'| \le
6n+4-(4\lceil\sqrt{n}\rceil +3) \le
6n-4\lceil\sqrt{n}\rceil +1.
$$
\end{proof}

Note that, for $5 \le n \le 25, 6n-4\lceil\sqrt{n}\rceil +1 \ge 6n-20$,
the upper bound established in \cite{HSV99} for any RVG. Hence Theorem~\ref{thm:bound}
says something stronger than that result only for $n\ge26$.

If $B_n$ denotes the maximum number of edges among all URVGs with $n$ vertices,
then Theorem~\ref{thm:bound} says that $B_n \le 6n-4\lceil\sqrt{n}\rceil +1$.
The next theorem uses  examples of URVGs with $n\ge64$ vertices and at least
$6n -12\bigl\lfloor\sqrt{n}\bigl\rfloor+6$ edges to  establish that
$B_n = 6n - \Theta(\sqrt{n})$.

\begin{thm}\label{thm:tightbound}
There is a URVG on $n$ vertices with at least 
$6n -12\bigl\lfloor\sqrt{n}\bigl\rfloor+6$ edges
for each  $n  \ge 64$.
\end{thm}

\begin{proof}

First consider the case when $n = k^2$ is a perfect square. Fig.~\ref{fig:examplelayout},
adapted from a figure in \cite{HSV99}, shows a layout of a URVG with 64 vertices.
The numbers on the squares indicate the degrees of the corresponding vertices, and
this example can be extended to give an analogous URVG with $k^2$ vertices for $k\ge8$.
For $k\ge8$, the resulting URVG on $n=k^2$ vertices has four vertices of degree 4, four
of degree 6, $4(k-3)$  of degree 7, four  of degree 10, $4(k-4)$ of degree 11, and
the rest of degree 12. Hence the degree sum of this graph is
$16+24+28(k-3)+40+44(k-4)+12(k^2-8k+16)=12(k-1)^2$, and therefore the graph
has $6(k-1)^2 =6n-12\sqrt{n}+6$ edges. Call this graph $G_n$.

\begin{figure}
  \centerline{\epsfig{file=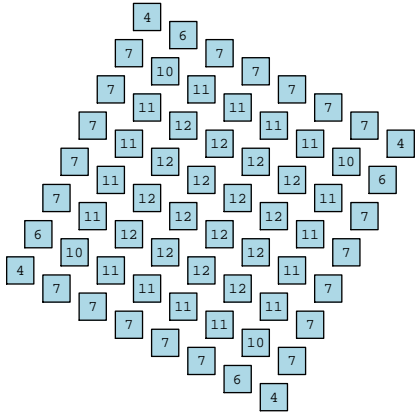, width=.45\textwidth
}}
\caption{A URV layout with $6n-12\sqrt{n}+6$ edges.}
\label{fig:examplelayout}
\end{figure}

Now suppose $k\ge8$. We show how to add additional unit squares to
an initial layout with $n=k^2$ squares to give examples satisfying the
condition for $n = k^2+1, k^2+2, \ldots, k^2+2k = (k+1)^2-1$. The numbered
squares in Fig.~\ref{fig:addingsquares} are being added to a layout with $n=64=8^2$,
but the method works for any $n=k^2$ with $k \ge 8$. We begin by adding the
northeasternmost square labeled 1, which adds one vertex and six new edges to the
corresponding graph. We then add, in sequence, the squares labeled $2, 3, 4, \ldots, k-4$,
completing a new row of $k-4$ squares along the upper side of the original layout.
Each vertex, when added, increases the edge count by 6. We repeat this process on
each of the other three sides of the original layout, adding a total of $4(k-4)$
squares, which is at least $2k$ squares, since $k\ge8$.

The function $f(n) = 6n -12\bigl\lfloor\sqrt{n}\bigl\rfloor+6$ increases
less quickly than the linear function $6n+6$, hence each unit increase of $n$
results in an increase of at most 6 for $f(n)$. Therefore, for $k\ge8, n = k^2$ and
$i = 0, \ldots, 2k-1$, we have by induction that
$$
f(n+i+1) \le f(n+i) + 6 \le |E(G_{n+i})|+6\le|E(G_{n+i+1})|
$$

\begin{figure}
  \centerline{\epsfig{file=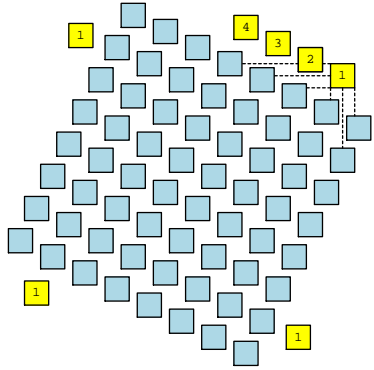, width=.4\textwidth}}
\caption{Adding squares to a layout with $k^2$ squares.}
\label{fig:addingsquares}
\end{figure}

This establishes the claim of the theorem.
\end{proof}

%\subsection{Edge bounds for bipartite URVGs}
\noindent \textbf{Edge bounds for bipartite URVGs} \\

In (\cite{DH94}, Corollary 10) it is shown that if a bipartite graph $G$ with $n\ge4$
vertices  is a subgraph of an RVG, then $G$ has at most $4n-12$ edges.
We use that result to bound the number of edges in a bipartite URVG.

\begin{thm}\label{thm:bipart}
For $n\ge 7$, let $G$ be a bipartite URVG with $n$ vertices. Then
$|E(G)| \le 4n-2\lceil\sqrt{n}\rceil+5$.
\end{thm}

\begin{proof}
Let $G$ be a bipartite graph with $n$ vertices and bipartition $\{A, B\}$,
and suppose that $L$ is a URV layout of $G$. As in Fig.~\ref{fig:enclosedlayout},
surround $L$ with the four rectangles $N, S, E, W$, creating an RVG layout that induces
a (non-bipartite) graph on $n+4$ vertices. 
As in the proof of Theorem~\ref{thm:bound} the number of edges with
exactly one endpoint in $\{N,S,E,W\}$ is at least $4\lceil\sqrt{n}\rceil -3$.
For each of the four rectangles, $\{N,S,E,W\}$, add it to $A$ if it sees 
more rectangles from $B$; otherwise add it to $B$.  
Call the enlarged bipartite sets $A'$ and $B'$.  
 
Define $G'$ to be the bipartite graph induced by the edges with one
endpoint in $A'$ and the other in $B'$.  All the edges of $G$ are edges
of $G'$ and there are at least $2\lceil\sqrt{n}\rceil -1$ 
bipartite edges with exactly one endpoint in $\{N,S,E,W\}$.
By a result of (\cite{DH94}, Corollary 10), the graph $G'$ has at most $4(n+4)-12=4n+4$ edges.

If we let $E'=E(G')\setminus E(G)$, then we have that
$$
|E(G)| \le 4n+4 - |E'| \le 4n+4-(2\lceil\sqrt{n}\rceil -1) =4n-2\lceil\sqrt{n}\rceil+5
$$
\end{proof}

Theorem~\ref{thm:tightbound} constructed examples of URVGs with
$6n - \Theta(\sqrt{n})$ edges. A construction similar to the one  in the proof of \
that theorem  gives examples, for $n\ge81$, of bipartite graphs with
$4n - \Theta(\sqrt{n})$ edges that are URVGs. %Fig.~\ref{fig:addingsquaresbi}
%shows the corresponding construction.

\begin{figure}
  \centerline{\epsfig{file=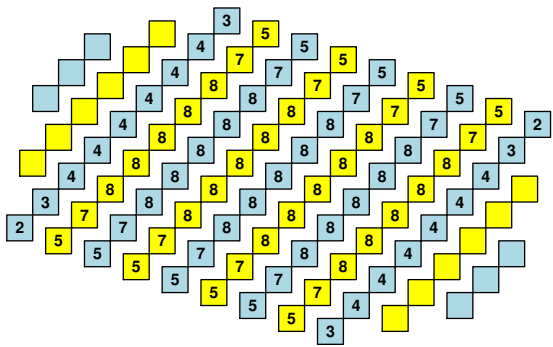, width=.6\textwidth}}
\caption{Adding squares to a weak bipartite layout with $k^2$ squares.}
\label{fig:addingsquaresbi}
\end{figure}

\begin{thm}\label{thm:tightbound2}
For every  $n  \ge 81$, there is a bipartite graph  with $n$ vertices
that is a  URVG and has  at least $4n -8\bigl\lceil\sqrt{n}\bigl\rceil+4$ edges.
\end{thm}

\begin{proof}
As in the proof of Theorem~\ref{thm:tightbound}, we do an initial construction
when $n=k^2$, and then we show how to handle the cases $n=k^2+1,\ldots,k^2+2k=(k+1)^2-1$.
The construction is illustrated with $k=9$ in Fig.~\ref{fig:addingsquaresbi}. 
The initial array of $k^2$ squares
(labeled with vertex degrees) has $4k^2-8k+4=4n-8\sqrt{n}+4$ edges. We then add the two
rows with $k-3$ unlabeled squares, each of which adds four new edges, and then the two
rows with $k-6$ unlabeled squares, adding four more edges each. We thus add a total
of $2(k-3)+2(k-6)=4k-18$ which is at least $2k$ if $k\ge9$, and the result follows
by the same argument as in the proof of Theorem~\ref{thm:tightbound}.
\end{proof}

\section{Open Problems}\label{sec:concl}
We conclude with two open problems motivated by the work in this paper.
\begin{enumerate}
  \item Theorem~\ref{thm:URVGtree} says that a tree is a URVG if and only
  if it is the union of two subdivided caterpillar forests, each with maximum
  degree 3. The same statement for general graphs is false, because, for example,
  there are URVGs with maximum degree greater than 6. However, we are not aware
  of any example of a graph $G$ that is the union of two subdivided caterpillar
  forests, each with maximum degree 3, that is \emph{not} a URVG. Note that
  Theorem~\ref{thm:linarb2} implies that any such example could not be constructed
  using two forests of paths. 
  %We conjecture that such examples do exist.
  \item While Theorem~\ref{thm:URVGtree} gives a complete characterization of
  URVG trees, it does not provide an efficient algorithm to decide whether an
  arbitrary tree is a URVG. Equivalently, no efficient algorithm is known for
  deciding if a tree is a union of two subdivided caterpillar forests, each
  with maximum degree 3. Peroche (\cite{Peroche82}, Theorem 4) has shown that deciding
  linear arboricity 2 is NP-complete. Shermer (\cite{Shermer96b}, Theorems 3.2 and 4.5) has shown that
  deciding caterpillar arboricity 2 is NP-complete, and also that deciding if a
  graph is an RVG is NP-complete. 
  %We conjecture that this problem is also NP-complete.
\end{enumerate}

\bibliography{bibliography}
\bibliographystyle{abbrv} % use bibliography style abbrv
\end{document}